\documentclass{article}
\usepackage{amsfonts}
\usepackage{amssymb}
\usepackage{graphicx}
\usepackage{amsmath}

\numberwithin{equation}{section} 
\newcommand{{\rr}}{{\mathbb R}}

\newenvironment{@abssec}[1]{
     \if@twocolumn
       \section*{#1}
     \else
       \vspace{.05in}\footnotesize
       \parindent .2in
         {\upshape\bfseries #1. }\ignorespaces
     \fi}
     {\if@twocolumn\else\par\vspace{.1in}\fi}
\newcommand\keywordsname{Key words}
\newcommand\AMSname{AMS subject classifications}

\begin{document}

\title{\textbf{The Set of Prime Numbers\thanks{%
Preprint n. $104$ deposited to Biblioteca Provinciale di Salerno
on 10 September 2007}}}
\author{Gerardo Iovane \\
Dipartimento di Ingegneria dell'Informazione e Matematica Applicata,\\
Universit\'{a} di Salerno\\
Via Ponte don Melillo, 84084 Fisciano (Sa), Italy,\\
iovane@diima.unisa.it\\
}
\date{\textbf{Submitted to Annals of Mathematics on }10.09.2007}
\maketitle

\begin{abstract}
In this work we show that the prime distribution is deterministic. Indeed
the set of prime numbers $P$\ can be expressed in terms of two subsets of $%
\mathbb{N}$ using three specific selection rules, acting on two sets of
prime candidates. The prime candidates are obtained in terms of the first
perfect number. The asymptotic behaviour is also considered.

We obtain for the first time an explicit relation for generating the full
set $P$ of prime numbers smaller than $n$ or equal to $n$\footnote{%
At the present the patent procedure is under definition. Consequently, each
a Computer Science procedure/algorithm, software solution or
firmware/hardware solution for generating prime candidates and selection
rules according to the results presented below for Public, Military/Defence
and Business aims, must be defined in agreement with the author of the
present paper. The solutions which are obtained from the relations shown in
what follows can be used for Public, Military/Defence and Business aims
exclusively in agreement with the author of the present paper.}.
\end{abstract}


\section{Introduction}

The sequence of prime numbers is of fundamental importance in many fields in
general, and in Mathematics in particular. During the last two centuries
many mathematicians have attempted to solve this problem using different
methods (see for example \cite{Bombieri}-\cite{Granville}).

The arithmetic of prime numbers has a crucial role in the present
Cryptography and Information Security. Indeed, many cryptosystems, such as
RSA, XTR, ECC (Elliptic Curve Cryptosystems), are based on our historical
ignorance about the inner nature of primality.

In addition, many papers are devoted to primes in Physics; a sort of
mathematical blueprint seems to guide the evolution of structures at
different length scales \cite{Elnaschie}, \cite{Iovane1}.

While in \cite{Iovane_primi} we built an explicating approach based on
dynamical processes and genetic algorithm, here we analyze the analytic
properties of the proposed sets of prime candidates. Then we consider the
selection rules to obtain two pure sets of primes, containing all prime
numbers with the exception of the first two.

A number of efficient algorithms have been known for a long time (for
details see \cite{Manidra}, \cite{Manidra2}, \cite{Miller}, \cite{Rabin},
\cite{Solovay}, \cite{Adleman}). The algorithms of Rabin, and Solovay and
Strassen are randomized. In addition the algorithm of Adleman et al.
requires (slightly) super-polynomial time, and the algorithm of Miller is in
P only under an unproved number-theoretic hypothesis \cite{Manidra2}. A
relevant contribution was given by Agrawal, Kayal and Saxena in 2004 \cite%
{Manidra}. Indeed, they proved that the problem is in P.\ We
analyze the asymptotic behaviour too and show that we can write a
deterministic second degree polynomial algorithm if we want
numerically verify our results to obtain primes.

In writing this work I was encouraged by a historical sentence by Fields
Medalist Prof. Bombieri, that is, more or less the following: 'When things
become too complex, sometime it makes sense to stop and ask: is my question
right?' Moreover, another encouragement comes from the well known sentence
in the \textit{Fine Hall} of Princeton: "Raffiniert ist der Herr Gott, aber
boshaft ist Er nicht". Then, the main idea used as starting point of this
work are:

\begin{itemize}
\item the use of a simple language to represent the objects of interest
(i.e. primes);

\item Nature manifests itself through beautiful symmetries and harmonies
based on primes in such a way that we can understand some interesting
results thanks to primality. \newline
\end{itemize}

The paper is organized as follows. In Sect.2 we give an useful partition of $%
\mathbb{N}$ in terms of the first perfect number, then we write the two
explicit maximum sets of prime candidates and using some specific selection
rules we reduce the previous sets to the two explicit maximum sets of prime
numbers; Sect.3 is devoted to study some asymptotic properties, while we
consider the algorithms and their computational complexity to verify the
results numerically in Sect.4. The conclusions are drawn in Sect.5.

\section{Prime Candidates and Prime Sets}

\subsection{The Sets of Prime Candidates}

It is well known a perfect number is defined as an integer which is the sum
of its proper positive divisors, that is, the sum of the positive divisors
not including the number. Equivalently, a perfect number is a number that is
half the sum of all of its positive divisors, or $\sigma (n)$ = 2 n. It is
very interesting to stress that the first perfect number is 6, because 1, 2
and 3 are its proper positive divisors and 1 + 2 + 3 = 6. For the reasons
which were expressed in \cite{Iovane_primi} let us write the set of natural
number $\mathbb{N}$ in terms of the first perfect number. Indeed, we
consider the following sets of positive integers

\begin{eqnarray}
ONE &=&\left\{ o_{k}=6k-5:k\in \mathbb{N}\right\} ,  \label{ONE} \\
TWO &=&\left\{ t_{k}=6k-4:k\in \mathbb{N}\right\} ,  \label{TWO} \\
THREE &=&\left\{ r_{k}=6k-3:k\in \mathbb{N}\right\} ,  \label{THREE} \\
FOUR &=&\left\{ f_{k}=6k-2:k\in \mathbb{N}\right\} ,  \label{FOUR} \\
A &=&\left\{ \alpha _{k}=6k-1:k\in \mathbb{N}\right\} ,  \label{A} \\
SIX &=&\left\{ s_{k}=6k:k\in \mathbb{N}\right\} ,  \label{SIX} \\
B &=&\left\{ \beta _{k}=6k+1:k\in \mathbb{N}\right\} ,  \label{B} \\
&&.....  \notag
\end{eqnarray}

\bigskip Then, we can write the set of natural numbers $\mathbb{N}$ as
follows

\begin{equation}
\mathbb{N}=\left\{ 1\right\} \cup TWO\cup THREE\cup FOUR\cup A\cup SIX\cup B
\label{N}
\end{equation}

Indeed, $ONE$ has the same elements of $B$ with the exception of the first
element, that is $1$. Moreover, also the sets which follows $B$ have got the
same elements of the previous sets ($ONE-B$), since their elements are just
obtained shifting the elements into the previous sets. The following table
shows the first values of the previous sets.

\begin{equation*}
\begin{tabular}{llllll}
$6k-4$ & $6k-3$ & $6k-2$ & $A$ & $6k$ & $B$ \\
2 & 3 & 4 & 5 & 6 & 7 \\
8 & 9 & 10 & 11 & 12 & 13 \\
14 & 15 & 16 & 17 & 18 & 19 \\
20 & 21 & 22 & 23 & 24 & 25 \\
26 & 27 & 28 & 29 & 30 & 31 \\
32 & 33 & 34 & 35 & 36 & 37 \\
... & ... & ... & ... & ... & ...%
\end{tabular}%
\end{equation*}

\bigskip \textbf{Lemma 1.} \textit{The set }$TWO$\textit{\ is made of
composite numbers with the exception of the first element, that is }$2$%
\textit{.}

\textbf{Proof.} The set $TWO$ contains elements of the form $6k-4$ for
construction. For each a $k\in \mathbb{N}$, $6k$ is an even number and the
same thing happens for $6k-4,$ since the difference between two even numbers
is an even number too. Then only the first element of $TWO$ can be prime and
this is the case due to the fact that the number $2$ is the first prime
number.\newline
$\square $ \newline

\textbf{Lemma 2.} \textit{The set }$THREE$\textit{\ is made of composite
numbers with the exception of the first element, that is }$3$\textit{.}

\textbf{Proof.} The set $THREE$ contains elements of the form $6k-3$ for
construction. For each a $k\in \mathbb{N}$, $6k$ is an even number and $6k-3$
is an odd multiple of $3$; that is each an element is odd and $6k-3=0$ (mod
3). Then only the first element of $THREE$ can be prime and this is the case
due to the fact that the number $3$ is the second prime number. \newline
$\square $

\bigskip \textbf{Lemma 3.} \textit{The sets }$FOUR$ and $SIX$\textit{\ are
made of \ composite numbers with no exception.}

\textbf{Proof.} The proof is trivial due to the fact that the elements of
the set $SIX$ are multiple of the first perfect number, that is an even
number (that is $6k=0$ (mod 6) ), while the elements of FOUR have the form $%
6k-2$ and the difference between two even numbers is an even number too.
\newline
$\square $

\bigskip \textbf{Theorem 1 (The Prime Candidates).} The set of the prime
candidates $\widetilde{P}$ can be written as

\begin{eqnarray}
\widetilde{P} &=&(\mathbb{N\cap }\left\{ 2\right\} )\cup (\mathbb{N\cap }%
\left\{ 3\right\} )\cup \left( \mathbb{N}\diagdown \left( TWO\cup THREE\cup
FOUR\cup SIX\right) \right) =  \notag \\
&=&\left\{ 2\right\} \cup \left\{ 3\right\} \cup A\cup B
\label{P_Candidates}
\end{eqnarray}

\textbf{Proof. }The proof follows from the Lemma 1-3 \textit{trivially.
Indeed note that for construction only the sets }$A$\textit{\ and }$B$%
\textit{\ among the others can be made of primes, since }$FOUR$\textit{\ and
}$SIX$\textit{\ are made of composite numbers, while }$TWO$\textit{\ and }$%
THREE$\textit{\ contain one prime only, that is the first element of each a
set.}\newline
$\square $ \newline

Unfortunately, the sets $A$ and $B$ also contain composite numbers such as
the positive integers which are multiple of $5,$ and so on. Consequently, we
introduce the selection rules for obtaining the full set of prime.

\subsection{Selection Rules and the Full Set of Prime Numbers}

\bigskip If we look at the sequence of prime candidates $x\geq 5$ coming
from $A$ and $B$ and compare them with the sequence of primes we see that,
with some exceptions, $A$ and $B$ alternate each others. For example, if we
consider $k=1,...,10$ we obtain

\begin{equation*}
\begin{array}{ccc}
& 6k-1 & 6k+1 \\
k=1 & 5 & 7 \\
k=2 & 11 & 13 \\
k=3 & 17 & 19 \\
k=4 & 23 & \mathbf{25} \\
k=5 & 29 & 31 \\
k=6 & \mathbf{35} & 37 \\
k=7 & 41 & 43 \\
k=8 & 47 & \mathbf{49} \\
k=9 & 53 & \mathbf{55} \\
k=10 & 59 & 61%
\end{array}%
\end{equation*}

where the numbers in bold face are composite. The total list of prime
numbers is in the following table. The second column of this table contain
the symbol "$-$" if the corresponding prime candidate comes from the set $A,$
and "$+$" if it comes from the set $B.$

\begin{equation*}
\begin{array}{cc}
5 & - \\
7 & + \\
11 & - \\
13 & + \\
17 & - \\
19 & + \\
23 & - \\
29 & - \\
31 & + \\
37 & + \\
41 & - \\
43 & + \\
47 & - \\
53 & - \\
59 & - \\
61 & +%
\end{array}%
\end{equation*}

It seems clear that the class transition, between two consecutive
primes, appears to be random. Indeed, in a first approximation it
can be modeled through a Random Process. But a deep analysis show
that the generation process is deterministic and can be realized
by a recursive combination of two sub-processes. \ In
\cite{Iovane_primi} the author showed that we must use a process
that produces a jump between the classes $A$ and $B$ (we can call
it \textit{zig-zag}) and a second process which switches off a
number
into a class for a fixed $k$ if it is composite (we can call it \textit{%
intermittence}). The intermittence can be made using the selection rules
presented below.

Let us introduce the following two subsets $A^{(-)}\subset A$ and $%
B^{(-)}\subset B$:

\begin{eqnarray}
A^{(-)} &=&\left\{ \alpha _{k_{ij}}\in A:\alpha
_{k_{ij}}=36ij-6i+6j-1,\forall k,i,j\in \mathbb{N}\right\} ,  \label{A-} \\
B^{(-)} &=&\left\{ \beta _{k_{ij}}\in B:\beta _{k_{ij}}=36ij-6i-6j+1\text{
and }\beta _{k_{ij}}=36ij+6i+6j+1\text{ }\forall k,i,j\in \mathbb{N}\right\}
.  \label{B-}
\end{eqnarray}

Then we can prove the following Lemmas.

\bigskip \textbf{Lemma 4}. \textit{The positive integer numbers }$\alpha
_{k_{ij}}\in A^{(-)}$\textit{\ are composite numbers . }

\textbf{Proof. }Multiplying $\alpha _{k}\in A$ to $\beta _{k}\in B$ we obtain

\begin{eqnarray*}
\alpha _{k_{ij}} &=&\alpha _{k}\beta _{k}=(6i+1)(6j-1)= \\
&=&36ij-6i+6j-1=6S_{ij}^{(1)}-1,
\end{eqnarray*}

where%
\begin{equation}
S_{ij}^{(1)}=(6ij-i+j)\text{ }\forall i,j\in \mathbb{N}.  \label{S1}
\end{equation}

Choosing $(6ij+i-j)=k$ we see that $(6i-1)(6j+1)=6k-1\in A.$ Consequently,
the product $\alpha _{k}\beta _{k}$ with $\alpha _{k}\in A$ and $\beta
_{k}\in B$ gives $\alpha _{k_{ij}}\in A^{\left( -\right) }$. This means that
the numbers $\alpha _{k_{ij}}\in A$ are composite numbers, since they can be
written as product of two numbers. $\square $ \newline

\textbf{Remarks 1.} The prime numbers $\alpha _{k}\in A$ must be different
from $\alpha _{k_{ij}}$.

\bigskip

\textbf{Lemma 5}. \textit{The positive integer numbers }$\beta _{k_{ij}}\in
B^{(-)}$\textit{\ are composite numbers.}

\textbf{Proof. }Multiplying two elements $\beta _{k},\beta _{k\ast }\in B$,
we obtain

\begin{eqnarray*}
\beta _{k_{ij}} &=&\beta _{k}\beta _{k\ast }=(6i+1)(6j+1)= \\
&=&36ij+6i+6j+1=6S_{ij}^{(2)}+1,
\end{eqnarray*}

where%
\begin{equation}
S_{ij}^{(2)}=(6ij+i+j)\text{ }\forall i,j\in \mathbb{N}.  \label{S2}
\end{equation}

Choosing $(6ij+i+j)=k$ we see that $(6i+1)(6j+1)=6k+1\in B.$ Consequently,
the product $\beta _{k}\beta _{k\ast }$ with $\beta _{k},\beta _{k\ast }\in
B $ gives $\beta _{k_{ij}}\in B^{\left( -\right) }$. This means that the
numbers $\beta _{k_{ij}}\in B$ are composite numbers, since they can be
written as product of two numbers.

Similarly, multiplying two elements $\alpha _{k},\alpha _{k\ast }\in A$, we
obtain

\begin{eqnarray*}
\beta _{k_{ij}} &=&\alpha _{k}\alpha _{k\ast }=(6i-1)(6j-1)= \\
&=&36ij-6i-6j+1=6S_{ij}^{(3)}+1,
\end{eqnarray*}

\bigskip where

\begin{equation}
S_{ij}^{(3)}=(6ij-i-j)\text{ }\forall i,j\in \mathbb{N}.  \label{S3}
\end{equation}%
Consequently, the product $\alpha _{k}\alpha _{k\ast }$ with $\alpha
_{k},\alpha _{k\ast }\in A$ gives $\beta _{k_{ij}}\in B$. This means that
the numbers $\beta _{k_{ij}}\in B^{\left( -\right) }$ are composite numbers
too, since they can be written as product of two numbers.

$\square $

\bigskip

\textbf{Remark 2.} The primes $\beta _{k}\in B$ must be different from $%
\beta _{k_{ij}}$.

\bigskip

\textbf{Remark 3}. \textit{The Lemmas 4-5 show that the numbers }$\alpha
_{k} $\textit{\ and }$\beta _{k}$\textit{, which are expressed as }$\alpha
_{k_{ij}}$\textit{\ and }$\beta _{k_{ij}}$ are composite numbers. Now let us
prove that they are the unique composite numbers in $A$ and $B$.\textit{\ }

\bigskip \textbf{Lemma 6}. \textit{The composite numbers }$\alpha
_{k_{ij}}\in A$\textit{\ can only be written in terms of the product }$%
\alpha _{k}\beta _{k\ast }.$

\textbf{Proof.} From the Lemma 5, it follows that \bigskip
\begin{equation*}
\alpha _{k}\alpha _{k\ast }=0\text{ (mod }6k+1\text{),}
\end{equation*}

that is $\alpha _{kk\ast }=\alpha _{k}\alpha _{k\ast }\in B$ \ $\forall
\alpha _{k},\alpha _{k\ast }\in A$. Moreover,
\begin{equation*}
t_{k}t_{k\ast }=f_{k}f_{k\ast }=s_{k}s_{k\ast }=t_{k}f_{k\ast
}=t_{k}s_{k\ast }=f_{k}s_{k\ast }=0\text{ (mod }2\text{).}
\end{equation*}%
that is the product of two elements of the sets $TWO,$ $FOUR$ and $SIX$ is
an even number.

In addition
\begin{equation*}
r_{k}r_{k\ast }=0\text{ (mod }3\text{) and it is an odd number,}
\end{equation*}%
that is $r_{k}r_{k\ast }\in THREE$. Furthermore,

\begin{equation*}
r_{k}t_{k\ast }=r_{k}f_{k\ast }=r_{k}s_{k\ast }=0\text{ (mod }2\text{),}
\end{equation*}%
that is the product of an element of the set $THREE$ to an element of the
sets $TWO$ (or $FOUR$ or $SIX$) is an even number. We also obtain that%
\begin{equation*}
\alpha _{k}t_{k\ast }=\alpha _{k}f_{k\ast }=\alpha _{k}s_{k\ast }=0\text{
(mod }2\text{),}
\end{equation*}%
that is the product of an $\alpha _{k}\in A$ to an element of the set $TWO$
(or $FOUR$ or $SIX$) is an even number, and also

\begin{equation*}
\alpha _{k}r_{k\ast }=0\text{ (mod }3\text{) and it is an odd number,}
\end{equation*}%
that is $\alpha _{k}r_{k\ast }\in THREE$. Moreover,\bigskip
\begin{equation*}
\beta _{k}\beta _{k\ast }=0\text{ (mod }6k+1\text{),}
\end{equation*}

that is $\beta _{kk\ast }=\beta _{k}\beta _{k\ast }\in B$ \ $\forall \beta
_{k},\beta _{k\ast }\in B$, and%
\begin{equation*}
\beta _{k}t_{k\ast }=\beta _{k}f_{k\ast }=\beta _{k}s_{k\ast }=0\text{ (mod }%
2\text{),}
\end{equation*}

\bigskip and%
\begin{equation*}
\beta _{k}r_{k\ast }=0\text{ (mod }3\text{) and it is an odd number.}
\end{equation*}

Consequently, we obtain that a composite element $\alpha _{k}\in A$\textit{\
can be written only in terms of the product }$\alpha _{k}\beta _{k\ast }$.

$\square $ \bigskip

\textbf{Lemma 7}. \textit{The composite numbers }$\beta _{k}\in B$ \textit{%
can only be written in terms of the products }$\alpha _{k}\alpha _{k\ast }$,
$\beta _{k}\beta _{k\ast }.$

\textbf{Proof.} The proof is obtained using the results of the Lemmas 4-6
and going head in the same way like in Lemma 6.

$\square $

\bigskip

\textbf{Theorem 2 (about the Selection Rules).} \textit{The natural numbers }%
$\alpha _{k}\in A$\textit{\ , }$\beta _{k}\in B$ \textit{are composite if
and only if }$\alpha _{k}\in A^{(-)}$\textit{\ , }$\beta _{k}\in B^{(-)}$%
\textit{. }

\textbf{Proof.} The proof of the theorem follows from the previous four
Lemmas trivially.

$\square $ \bigskip

From (\ref{N}), and the theorem 2 we obtain the following theorem.

\bigskip

\textbf{Theorem 3. (The Full Set of Primes).}\textit{\ The full set of
primes has got the following minimum explicit representation}%
\begin{equation}
P=\left\{ 2\right\} \cup \left\{ 3\right\} \cup A^{\prime }\cup B^{\prime },
\label{P}
\end{equation}

\textit{where}%
\begin{align}
A^{\prime }& =\left\{ \alpha _{k}\in \mathbb{N}:\alpha _{k}=6k-1\text{ and }%
k\neq 6ij-i+j,\forall k,i,j\in \mathbb{N}\right\} ,  \label{A'} \\
B^{\prime }& =\left\{ \beta _{k}\in \mathbb{N}:\beta _{k}=6k+1\text{ and }%
k\neq 6ij+i+j\text{ or }k\neq 6ij-i-j,\forall k,i,j\in \mathbb{N}\right\} .
\label{B'}
\end{align}

\bigskip \textbf{Remark 4.} \ We can also note that%
\begin{equation}
A^{\prime }=A\diagdown A^{(-)}\text{ and }B^{\prime }=B\diagdown B^{(-)}.
\label{R4.1}
\end{equation}

then%
\begin{equation}
P=\left\{ 2\right\} \cup \left\{ 3\right\} \cup (A\diagdown A^{(-)})\cup
(B\diagdown B^{(-)}).  \label{R4.2}
\end{equation}

\section{Some Asymptotic Consequences}

Let us consider some asymptotic behaviours. \newline

\textbf{Theorem 4 (Equipartition of the Cuts).} If we call $\#_{A}(s)$ the
first $s$\ elements of $A$ and $\#_{B}(s)$ the first $s$\ elements of $B$,
then $\#_{A^{(-)}}(s)=s^{2}$ and $\#_{B^{(-)}}(s)=s^{2}+s$. Moreover, for $s$
$\rightarrow $ $\infty $ \ we have $\#_{A^{(-)}}(s)\approx \#_{B^{(-)}}(s).$

\textbf{Proof.} Considering that each an element of $A^{(-)}$ is given
multiplying an element of $A$ to an element of $B$ we obtain $%
\#_{A^{(-)}}(s)=s^{2}$ trivially. Considering that each an element of $%
B^{(-)}$ is given multiplying two elements of $A$ or two elements of $B$,
then

$\#_{B^{(-)}}(s)=\left(
\begin{array}{c}
s+2-1 \\
2%
\end{array}%
\right) +\left(
\begin{array}{c}
s+2-1 \\
2%
\end{array}%
\right) =$ $s^{2}+s.$

Moreover,

\begin{equation}
\lim_{k\rightarrow \infty }\frac{\#_{B^{(-)}}(s)}{\#_{A^{(-)}}(s)}%
=\lim_{k\rightarrow \infty }\frac{s^{2}+s}{s^{2}}=1\text{.}
\end{equation}
$\square $

\bigskip \textbf{Theorem 5 (Equipower of Primes into the Sets of Primes).}
If we call $\#_{A^{\prime }}(k)$ the number of the first $k$ primes into $%
A^{\prime },$ and $\#_{B^{\prime }}(k)$ the number of the first $k$ primes
into $B^{\prime }$, then

\begin{equation*}
\lim_{k\rightarrow \infty }\frac{\#_{B^{\prime }}(k)}{\#_{A^{\prime }}(k)}=l,%
\text{ with }l<\infty
\end{equation*}

\textbf{Proof.} For a fixed $k$ we obtain $k$ elements into $A$, that is $%
\#_{A}(k)=k;$ while to evaluate the number of selectors we must estimate the
maximum indexes $i$ and $j$. Indeed from (\ref{S1}) we obtain%
\begin{equation}
i_{\text{max}}^{(S1)}=\left\lfloor \frac{k-1}{5}\right\rfloor ,
\label{i_max_s1}
\end{equation}

corresponding to $j=0$ and
\begin{equation}
j_{\text{max}}^{(S1)}=\left\lfloor \frac{k+1}{7}\right\rfloor ,
\label{j_max_s1}
\end{equation}

corresponding to $i=0$. Consequently the number of selectors smaller then $k$
or equal to $k$, with respect to (\ref{S1}), is
\begin{equation}
\maltese ^{(S1)}(k)=i_{\text{max}}^{(S1)}j_{\text{max}}^{(S1)}=\left\lfloor
\frac{k-1}{5}\right\rfloor \left\lfloor \frac{k+1}{7}\right\rfloor .
\label{num_S1}
\end{equation}

Then the number of primes into $A$ for a fixed $k$ will be%
\begin{eqnarray}
\#_{A^{\prime }}(k) &=&\#_{A}(k)-\maltese ^{(S1)}(k)=  \notag \\
&=&k-\left\lfloor \frac{k-1}{5}\right\rfloor \left\lfloor \frac{k+1}{7}%
\right\rfloor  \label{PHI_A'}
\end{eqnarray}

For the same $k$ we obtain $k$ elements into $B$ too, that is
$\#_{B}(k)=k;$ while to evaluate the number of selectors we must
estimate the maximum indexes $i$ and $j$ coming from (\ref{S2})
and (\ref{S3}). Indeed for the selectors of $B$, we stress that
when $i$ and $j$ run from $1$ to their max value we obtain a
symmetric matrix of selectors. Then arbitrarily for one of the two
indexes we must choose the maximum with $\left\lceil \cdot
\right\rceil $, instead of $\left\lfloor \cdot \right\rfloor $. In
other words from (\ref{S2}), we have

\begin{equation}
i_{\text{max}}^{(S2)}=\left\lfloor \frac{k-1}{7}\right\rfloor ,
\label{i_max_s2_1}
\end{equation}

and
\begin{equation}
j_{\text{max}}^{(S2)}=\left\lceil \frac{k-1}{7}\right\rceil .
\label{imax_s2_1}
\end{equation}

Consequently the number of selectors smaller then $k$ or equal to $k $, with
respect to (\ref{S2}), is
\begin{equation}
\maltese ^{(S2)}(k)=i_{\text{max}}^{(S2)}j_{\text{max}}^{(S2)}=\left\lfloor
\frac{k-1}{7}\right\rfloor \left\lceil \frac{k-1}{7}\right\rceil .
\label{num_S2}
\end{equation}%
Moreover from (\ref{S3}), we obtain

\begin{equation}
i_{\text{max}}^{(S3)}=\left\lfloor \frac{k+1}{5}\right\rfloor ,
\label{i_max_s2_2}
\end{equation}

and%
\begin{equation*}
j_{\text{max}}^{(S3)}=\left\lceil \frac{k+1}{5}\right\rceil .
\end{equation*}

Then the number of selectors smaller then $k$ or equal to $k$, with respect
to (\ref{S3}), is%
\begin{equation}
\maltese ^{(S3)}(k)=i_{\text{max}}^{(S3)}j_{\text{max}}^{(S3)}=\left\lfloor
\frac{k+1}{5}\right\rfloor \left\lceil \frac{k+1}{5}\right\rceil .
\label{num_S3}
\end{equation}

Consequently the number of primes into $B$ for a fixed $k$ will be
\begin{eqnarray}
\#_{B^{\prime }}(k) &=&\#_{B}(k)-\maltese ^{(S2)}(k)-\maltese ^{(S3)}(k)=
\notag \\
&=&k-\left\lfloor \frac{k-1}{7}\right\rfloor \left\lceil \frac{k-1}{7}%
\right\rceil -\left\lfloor \frac{k+1}{5}\right\rfloor \left\lceil \frac{k+1}{%
5}\right\rceil .  \label{PHI_B'}
\end{eqnarray}

Without losing generality for our purpose, let us approximate
\begin{equation}
\left\lfloor \frac{k\pm 1}{\delta }\right\rfloor =\left\lceil \frac{k\pm 1}{%
\delta }\right\rceil \approx \left( \frac{k\pm 1}{\delta }\right) \text{
with }\delta =5,7;  \label{APPROX}
\end{equation}

\begin{eqnarray}
\lim_{k\rightarrow \infty }\frac{\#_{B^{\prime }}(k)}{\#_{A^{\prime }}(k)}
&=&\lim_{k\rightarrow \infty }\frac{k-\left\lfloor \frac{k-1}{7}%
\right\rfloor \left\lceil \frac{k-1}{7}\right\rceil -\left\lfloor \frac{k+1}{%
5}\right\rfloor \left\lceil \frac{k+1}{5}\right\rceil }{k-\left\lfloor \frac{%
k-1}{5}\right\rfloor \left\lfloor \frac{k+1}{7}\right\rfloor }\approx
\notag \\
&\approx &\lim_{k\rightarrow \infty }\frac{\mu k^{2}+\lambda k+\omega }{\rho
k^{2}+\vartheta k+\sigma }=\mu /\rho =l,  \label{ASYMP}
\end{eqnarray}

with\ $\mu ,\lambda ,\omega ,\rho ,\vartheta ,\sigma $-constants; in
particular in the present case $\mu =74$ and $\rho =1;$ this means $%
\#_{A^{\prime }}(k)\approx \frac{1}{\mu }\#_{B^{\prime }}(k)~$for $%
k\rightarrow \infty .$

$\square $

\bigskip

\textbf{Remark 5.} \textit{Of course the approximation (\ref{APPROX}) is not
acceptable for evaluating $\pi (n)$. Indeed in this case we could solve the
inverse problem: the quantities $\pi (n)$ and $\#(n)$ could be considered as
known quantities for estimating $\maltese (n)$ (for details see \cite%
{Iovane_primi})}.

\section{Algorithms and their Computational Complexity}

If we use the definition of prime numbers we have immediately a way of
determining if a number $n$ is prime. Indeed, try dividing $n$ by every
number $m\leq \sqrt{n}$ then if any $m$ divides $n$ then the last one is
composite, otherwise it is prime. As it is well known this test is
inefficient, since it takes $\Omega (\sqrt{n})$ steps to determine if $n$ is
prime. As anticipated and shown in \cite{Manidra} an unconditional
deterministic polynomial-time algorithm, for determining whether an input is
prime or composite, can be obtained. In what follows we show that our
not-optimized algorithm has a computational complexity $C(n)\in O(n^{2})$ if
we want verify whether $n$ is prime or composite, while a second algorithm
to verify the same things according to the sets $A^{\prime }$ and $B^{\prime
}$ has a computational complexity $C(n)\in O(1)$.

It is trivial to prove that following the transformation $r=i-1,$ $s=j-1$,
the previous three selection rules (\ref{S1}), (\ref{S2}), (\ref{S3}) can be
written respectively%
\begin{eqnarray}
S_{rs}^{(1)} &=&6+7s+(5+6s)r\text{ \ }\forall r,s\in \mathbb{N}_{0},
\label{S1_1} \\
S_{rs}^{(2)} &=&8+7s+(7+6s)r\text{ \ }\forall r,s\in \mathbb{N}_{0},
\label{S2_1} \\
S_{rs}^{(3)} &=&4+5s+(5+6s)r\text{\ }\forall r,s\in \mathbb{N}_{0}.
\label{S3_1}
\end{eqnarray}

The not-optimized algorithm to verify whether $n$ is prime can be written as
follows.

\begin{equation*}
\begin{tabular}{|l|}
\hline
$%
\begin{array}{l}
\text{Input: integer }n>3 \\
\text{Step 1. If \ }(n+1=0\text{ (}mod\text{ }6\text{) })\text{ \ }k=(n+1)/6,%
\text{ }r_{\text{max}}=\left\lfloor \frac{k-6}{5}\right\rfloor ,\text{ }s_{%
\text{max}}=\left\lfloor \frac{k-6}{7}\right\rfloor  \\
\text{ \ \ \ \ \ \ \ \ \ \ \ \ \ \ \ \ For }s=0\text{ to }s_{\text{max}}%
\text{ do} \\
\text{ \ \ \ \ \ \ \ \ \ \ \ \ \ \ \ \ \ \ \ \ \ \ For }r=0\text{ to }r_{%
\text{max}}\text{ do} \\
\text{ \ \ \ \ \ \ \ \ \ \ \ \ \ \ \ \ \ \ \ \ \ \ \ \ \ \ \ \ }%
k(r,s)=6+7s+(5+6s)r \\
\text{ \ \ \ \ \ \ \ \ \ \ \ \ \ \ \ \ \ \ \ \ \ \ \ \ \ \ \ \ If }(k(r,s)=k)%
\text{ output COMPOSITE} \\
\text{ \ \ \ \ \ \ \ \ \ \ \ \ \ \ \ \ \ \ \ \ \ \ \ \ \ \ \ \ Else output
PRIME} \\
\text{Step 2. If \ }(n-1=0\text{ (}mod\text{ }6\text{) })\text{ \ }k=(n-1)/6,%
\text{ }l_{\text{max}}=\left\lfloor \frac{k-4}{5}\right\rfloor ,\text{ }m_{%
\text{max}}=l_{\text{max}}+1 \\
\text{ \ \ \ \ \ \ \ \ \ \ \ \ \ \ \ \ For }m=0\text{ to }m_{\text{max}}%
\text{ do} \\
\text{ \ \ \ \ \ \ \ \ \ \ \ \ \ \ \ \ \ \ \ \ \ For }l=m\text{ to }l_{\text{%
max}}\text{ do} \\
\text{\ \ \ \ \ \ \ \ \ \ \ \ \ \ \ \ \ \ \ \ \ \ \ \ \ \ \ }%
k(l,m)=4+5m+(5+6m)l \\
\text{ \ \ \ \ \ \ \ \ \ \ \ \ \ \ \ \ \ \ \ \ \ \ \ \ \ \ If }(k(l,m)=k)%
\text{ output COMPOSITE} \\
\text{ \ \ \ \ \ \ \ \ \ \ \ \ \ \ \ \ \ Else }c_{\text{max}}=\left\lfloor
\frac{k-8}{7}\right\rfloor ,\text{ }d_{\text{max}}=d_{\text{max}}+1 \\
\text{ \ \ \ \ \ \ \ \ \ \ \ \ \ \ \ \ \ For }d=0\text{ to }d_{\text{max}}%
\text{ do} \\
\text{ \ \ \ \ \ \ \ \ \ \ \ \ \ \ \ \ \ \ \ \ \ \ \ For }c=d\text{ to }c_{%
\text{max}}\text{ do} \\
\text{ \ \ \ \ \ \ \ \ \ \ \ \ \ \ \ \ \ \ \ \ \ \ \ \ \ \ \ \ }%
k(c,d)=8+7d+(7+6d)c \\
\text{ \ \ \ \ \ \ \ \ \ \ \ \ \ \ \ \ \ \ \ \ \ \ \ \ \ \ \ \ If }(k(c,d)=k)%
\text{ output COMPOSITE} \\
\text{ \ \ \ \ \ \ \ \ \ \ \ \ \ \ \ \ \ \ \ \ \ \ \ \ \ \ \ \ Else output
PRIME} \\
\text{Step 3 Else output COMPOSITE}%
\end{array}%
$ \\ \hline
\end{tabular}%
\end{equation*}

We are taking into account that the matrices of the selectors $k(l,m)$ and $%
k(c,d)$ are symmetric.

\bigskip

\textbf{Theorem 6.} \textit{The algorithm above returns PRIME if and only if
}$n$\textit{\ is prime.}

\textbf{Proof.} If $n\notin A$ or $n\notin B$ $Step$ $3$, and so the
algorithm, returns COMPOSITE. Moreover if $n\in A$ and $k(r,s)=k$ then $Step$
$1$ returns COMPOSITE. For each a $n\in B$ and $k(l,m)=k$ (or $k(c,d)=k$) $%
Step$ $2$ returns COMPOSITE. Consequently, the algorithm returns PRIME if
only $n\in A^{\prime }$ or $n\in B^{\prime },$ that is if only $n$ is prime.

Vice versa if $n$ is prime then $n\in A^{\prime }$ or $n\in B^{\prime }$.

If $n\in A^{\prime }$ then $n\in A$ since $A^{\prime }\subset A$. If $n\in A$
the algorithm returns COMPOSITE for each a $n$ such that $k(r,s)=k$, but if $%
n$ is prime, then $n\notin A^{(-)}.$ Then the output of $Step$ $1$ must be
PRIME.

Similarly if $n\in B^{\prime }$ then $n\in B$ since $B^{\prime }\subset B$.
If $n\in B$ the algorithm returns COMPOSITE for each a $n$ such that $%
k(l,m)=k$ or $k(c,d)=k$, but if $n$ is prime, then $n\notin B^{(-)}.$ Then
the output of $Step$ $2$ must be PRIME.

$\square $ \bigskip

By looking at the previous algorithm we can evaluate $O(n)$ for recognizing
if $n$ is prime or composite.

This algorithm is not optimized, since our aim is just to show that the
computational complexity is a second degree deterministic polynom. \bigskip

\bigskip

\textbf{Theorem 7.} \textit{The computational complexity for recognizing if }%
$n$\textit{\ is prime or composite is}

\begin{equation}
C(n)\in O(n^{2})
\end{equation}

\textbf{Proof.} When $n$ $\in A$ is prime, the number of operations can be
evaluated as follows. The Step 1 costs two basic operations (that is a sum
and a division), to verify $n$ $\in A.$ In addition we have six basic
operations to estimate the indexes maximum values and $k$, $6\times r_{\text{%
max}}\times s_{\text{max}}$ operations to evaluate $k(r,s)$ and other $r_{%
\text{max}}\times s_{\text{max}}$ operations to compare $k(r,s)$ with $k.$
Consequently, $O_{\text{min}}(r_{\text{max}},s_{\text{max}})=8+7\times r_{%
\text{max}}\times s_{\text{max}}.$Considering that $r_{\text{max}%
}=\left\lfloor \frac{k-6}{5}\right\rfloor =\left\lfloor \frac{n-35}{5}%
\right\rfloor $ and $s_{\text{max}}=\left\lfloor \frac{k-6}{7}\right\rfloor
=\left\lfloor \frac{n-35}{42}\right\rfloor $ we obtain

\begin{equation}
C(n)\in O(n^{2})
\end{equation}

Similarly when $n$ $\in B$ is prime, the Step 2 requires two basic
operations to verify $n$ $\in B.$ In addition we have five basic operations
to estimate the indexes maximum values and $k$, $7\times r_{\text{max}%
}\times (r_{\text{max}}+1)$ operations to evaluate $k(r,s)$ and to compare $%
k(r,s)$ with $k.$ Consequently, $O_{\text{max}}(r,s)=7+7\times r_{\text{max}%
}\times (r_{\text{max}}+1).$ Considering that $r_{\text{max}}=\left\lfloor
\frac{k-4}{5}\right\rfloor =\left\lfloor \frac{n-25}{30}\right\rfloor $ or $%
r_{\text{max}}=\left\lfloor \frac{k-8}{7}\right\rfloor =\left\lfloor \frac{%
n-49}{42}\right\rfloor ,$ it follows

\begin{equation}
C(n)\in O(n^{2})
\end{equation}

$\square $

\bigskip

\textbf{Remark 6.} \textit{If }$n$\textit{\ is composite, then }$C(n)\in
\Omega (1)$\textit{.}

\bigskip

\textbf{Remark 9.} \textit{Theorem 7 is useful just to verify that }$C(n)\in
O(n^{2}).$

\bigskip

To evaluate if $n$ is prime or composite we use the following algorithm
based on the knowledge of $A^{\prime }$ and $B^{\prime }.$

\bigskip
\begin{equation*}
\begin{tabular}{|l|}
\hline
$%
\begin{array}{l}
\text{Input: integer }n>3 \\
\text{Step 1. If \ }(n+1=0\text{ (}mod\text{ }6\text{) })\text{ \ }k=(n+1)/6
\\
\text{ \ \ \ \ \ \ \ \ \ \ \ \ \ \ \ \ read }k(r,s) \\
\text{ \ \ \ \ \ \ \ \ \ \ \ \ \ \ \ \ If }(k(r,s)=k)\text{ output COMPOSITE}
\\
\text{ \ \ \ \ \ \ \ \ \ \ \ \ \ \ \ \ Else output PRIME} \\
\text{Step 2. If \ }(n-1=0\text{ (}mod\text{ }6\text{) })\text{ \ }k=(n-1)/6
\\
\text{ \ \ \ \ \ \ \ \ \ \ \ \ \ \ \ \ read }k(l,m) \\
\text{ \ \ \ \ \ \ \ \ \ \ \ \ \ \ \ \ If }(k(l,m)=k)\text{ output COMPOSITE}
\\
\text{ \ \ \ \ \ \ \ \ \ \ \ \ \ \ \ \ \ Else read }k(c,d) \\
\text{ \ \ \ \ \ \ \ \ \ \ \ \ \ \ \ \ \ If }(k(c,d)=k)\text{ output
COMPOSITE} \\
\text{ \ \ \ \ \ \ \ \ \ \ \ \ \ \ \ \ \ Else output PRIME} \\
\text{Step 3.\ \ Else output COMPOSITE }%
\end{array}%
$ \\ \hline
\end{tabular}%
\end{equation*}

\bigskip \textbf{Theorem 8 (Primality Test based on Pre-computed Selection
Rules). } \textit{The computational complexity for recognizing if }$n$%
\textit{\ is prime or composite is}

\begin{equation}
C(n)\in O(1)
\end{equation}

\textbf{Proof.} The proof is trivial using the rule of the sum for the
computational complexity and assuming that the arithmetic operations, the
function print and the function read have a computational complexity $%
C(n)\in O(1)$.

$\square $

\bigskip

\textbf{Remark 10.} \textit{Clearly we have to take into account the
computational complexity to generate the pre-computed selectors }$k$\textit{%
. Trivially considering that we have two nested for loops, it is in }$%
O(n^{2}).$\textit{\ We must stress this evaluation is off-line. In other
words, it is not in the procedure of testing, but it is into the procedure
of generation.}

\section{Conclusion}

In this work we defined two maximum explicit sets of prime candidates. Then
we showed that they can be reduced to two maximum explicit sets of primes.
In conclusion the sets obtained thanks to $6k\pm 1$ $k\in \mathbb{N},$ with
their selection rules , that is

\textit{where}%
\begin{align}
A^{\prime }& =\left\{ \alpha _{k}\in \mathbb{N}:\alpha _{k}=6k-1\text{ and }%
k\neq 6ij-i+j,\forall k,i,j\in \mathbb{N}\right\} ,  \\
B^{\prime }& =\left\{ \beta _{k}\in \mathbb{N}:\beta _{k}=6k+1\text{ and }%
k\neq 6ij+i+j\text{ or }k\neq 6ij-i-j,\forall k,i,j\in
\mathbb{N}\right\} .
\end{align}

together with $\left\{ 2\right\} $ and $\left\{ 3\right\} $ give us the full
set of prime numbers (in its reduced form), that is
\begin{equation*}
P=\left\{ 2\right\} \cup \left\{ 3\right\} \cup A^{\prime }\cup B^{\prime }.
\end{equation*}

We also show the selection rules have the same weight acting on the two sets
of candidates and the sets of primes have the same power.

Moreover basing on the results shown above, the algorithm is a second degree
deterministic polynomial procedure.

Thanks to the discovery \ of the sets $A^{\prime }$ and $B^{\prime },$ we
have obtained, for the first time an explicit expression of the full set $P$
of prime numbers smaller than $n$ or equal to $n$, which are generated with
a specific rule and without the use of some test.

\bigskip \textbf{Acknowledgements}

The author wishes to thank prof.Saverio Salerno, who stimulated him to
investigate primes.

\end{document}